\documentclass{amsart}

\newtheorem{thm}{Theorem}[]

\theoremstyle{definition}

\newtheorem{example}[thm]{Example}

\theoremstyle{remark}

\numberwithin{equation}{section}

\newcommand{\R}{\mathbb{R}} 

\begin{document}

\title{Dual Dictionaries in Linear Programming}

\author{Patrick T. Perkins}
\address{Department of Mathematics, University of Washington, 
Seattle, WA 98195}
\email{pperkins@uw.edu}

\author{Xiang Gao}
\address{Department of Mathematics, University of Washington, 
Seattle, WA 98195}
\email{seangao@uw.edu}

\date{}

\begin{abstract}
In order to use the Dual Simplex Method, one needs to prove a certain bijection between the dictionaries associated with the primal problem and those associated with its dual.  We give a short conceptual proof of why this bijection exists.
\end{abstract}

 \maketitle

\section{Introduction}

Chv\'{a}tal \cite{chv} introduces the notion of a {\it dictionary} associated to a Linear Programming problem (LP).  In order to use the Dual Simplex Method, one needs to prove a certain bijection between the dictionaries associated with the primal problem and those associated with its dual. Chv\'{a}tal leaves the proof as an exercise, involving a long computation. Vanderbei \cite{vdb} gives a short and elegant proof.  Our contribution is a short proof that, we feel, gives a clear conceptual reason for why this beautiful bijection exists.

First, we set up some notation we will use throughout the paper. Consider a general LP problem
\begin{equation}
\begin{aligned}
& {\text{max}}
& & z=\mathbf{c}^T \mathbf{x} \\
& \text{s.t.} & &  A_0\mathbf{x} \leq \mathbf{b} \\
& & &  \mathbf{x} \geq 0\\
\end{aligned}
\end{equation}
The dual problem is
\begin{equation}
\begin{aligned}
& {\text{max}}
& & -w=-\mathbf{b}^T \mathbf{y} \\
& \text{s.t.} & &  -A_0^T\mathbf{y} \leq -\mathbf{c} \\
& & &  \mathbf{y} \geq 0\\
\end{aligned}
\end{equation}
Here ${\bf x}\in\R^n$, ${\bf y}\in\R^m$ and $A_0$ is an $m\times n$ matrix. 
But we immediately introduce slack variables and, for the rest of the paper, take ${\bf x},{\bf y}\in\R^{m+n}$. Write $A=[\,A_0\, I\,]$ for the larger matrix with an $m\times m$ identity matrix appended to $A_0$. As usual, $x_1,\ldots,x_n$ are the decision variables for the primal problem and $x_{n+1},\ldots,x_{m+n}$ are its slack variables. But, following Chv\'{a}tal, we use $y_{n+1},\ldots,y_{n+m}$ as the decision variables for the dual problem and $y_1,\ldots,y_n$ for its slacks. This makes the bijection easier to see.

\begin{example}
{\rm If the initial dictionary for a primal problem is}
    \[
    \begin{aligned}
  x_4 &= 18 - 4x_1 - 2x_2 + 2x_3\\
   x_5 &= -3 + x_1 + x_2 + 2x_3\\
   z &= 8x_1 + 11x_2 -10x_3 \\
    \end{aligned}
    \]
        {\rm then the initial dictionary for the dual problem is}
    \[
    \begin{aligned}
    y_1 &= -8 + 4y_4 - y_5\\
    y_2 &= -11 + 2y_4 - y_5\\
    y_3 &= 10 - 2y_4 - 2y_5\\
    -w &= -18y_4 + 3y_5\\
    \end{aligned}
    \]
{\rm Pivoting once in the primal, letting $x_1$ enter the basis and $x_5$ leave, gives}
    \[
    \begin{aligned}
    x_4 &= 6 - 4x_5 + 2x_2 + 10x_3\\
    x_1 &= 3 + x_5 - x_2 - 2x_3\\
    z &= 24 + 8x_5 + 3x_2 -26x_3 \\
    \end{aligned}
    \]
    {\rm The corresponding pivot in the dual lets $y_5$ enter and $y_1$ leave.}
    \[
    \begin{aligned}
    y_5 &= -8 + 4y_4 - y_1\\
    y_2 &= -3 - 2y_4 + y_1\\
    y_3 &= 26 - 10y_4 + 2y_1\\
    -w &= -24 - 6y_4 - 3y_1\\
    \end{aligned}
    \]
{\rm Note that each dictionary for the dual LP is, in some sense, the negative transpose of the corresponding dictionary for the primal.}
\end{example}

To be more precise, let $B\cup N$ be an ordered partition of $\{1,2,...,m+n\}$ such that $|B|=m$ and the columns of $A=[A_0 \,I\,]$ indexed by $B$ are linearly independant.  Let $\mathbf{x}_B$ be the vector of variables indexed by $B$, and similarly for $\mathbf{x}_N$. Then the dictionary of the primal LP associated to this partition is of the form
\begin{equation}
\begin{aligned}
    \mathbf{x}_B &= \mathbf{p} - Q \,\mathbf{x}_N\\
    z &= z^*+\mathbf{q}^T \mathbf{x}_N
\end{aligned}
\end{equation}
where $Q$ is an $m\times n$ matrix, $\mathbf{p}\in\R^m$, $\mathbf{q}\in\R^n$ and $z^*\in\R$.

Given this set up, we will prove that 
\begin{equation}
\begin{aligned}
    \mathbf{y}_N &= -\mathbf{q} + Q^T \mathbf{y}_B\\
    -w &= -z^*-\mathbf{p}^T \mathbf{y}_B
\end{aligned}
\end{equation}
is a dictionary for the dual LP. This means that every solution to (1.4) is a solution to the initial dual dictionary, and vice versa.


\section{Orthogonal Subspaces}

We first recast our pair of LPs in terms of orthogonal subspaces.  This formulation is well known, we first encountered it in Todd \cite{todd}.  Add two new variables, $x_0$ and $x_{m+n+1}$, and set $\overline{\bf{x}}=[x_0,x_1, \ldots x_{m+n},x_{m+n+1}]^T\in\R^{n+m+2}$.  Define the matrix $R$ by
\[
R=\begin{bmatrix}
{\bf 0}^T & A_0 & I & -\mathbf{b}\\
1 & -\mathbf{c}^T & {\bf 0} & 0
\end{bmatrix} 
\]
Then the primal LP can be formulated
\begin{equation}
\begin{aligned}
& {\text{max}}
& & x_0 \\
& \text{s.t.} & &    x_1, \ldots x_{m+n} \geq 0\\
&&&x_{m+n+1}=1\\
&&&\overline{\bf x} \in\ker(R)
\end{aligned}\\
\end{equation}

Now consider the row space of $R$. Let $[\mathbf{u}^T,u_0]^T \in \mathbb{R}^{m+1}$. Every vector in the row space is of the form
\begin{equation}
[\mathbf{u}^T , u_0]
\begin{bmatrix}
{\bf 0}^T & A_0 & I & -\mathbf{b}\\
1 & -\mathbf{c}^T & {\bf 0} & 0
\end{bmatrix} =
\begin{bmatrix}
u_0 & \mathbf{u}^T A_0-u_0\mathbf{c}^T & \mathbf{u}^T & -\mathbf{u}^T \mathbf{b}
\end{bmatrix}
\end{equation}
Let $\overline{\bf{y}}=[y_0,y_1,\ldots,y_{m+n},y_{m+n+1}]\in\R^{m+n+2}$. Then we can reformulate the dual LP as
\begin{equation}
\begin{aligned}
& {\text{max}}
& & y_{m+n+1} \\
& \text{s.t.} & &    y_1, \ldots y_{m+n} \geq 0\\
&&&y_0=1\\
&&&\overline{\bf{y}}\in {\rm row\ space}(R)\\
\end{aligned}\\
\end{equation}
Thus $\R^{m+n+2}$ splits into two orthogonal subspaces, one associated with the primal LP and one with the dual.  Note that this naturally makes $y_n,\ldots,y_{m+n}$ the decision variables for the dual LP.
\section{The Proof}
Let ${\bf y}=[y_1,\ldots,y_{m+n}]^T\in\R^{m+n}$. Then $[{\bf y},w]^T$ is a solution to (1.4) if and only if 
\[
[1,{\bf y}_N^T,{\bf y}_B^T,w]\in {\rm row\ space}\left( \begin{bmatrix}
{\bf 0}^T & Q & I & -\mathbf{p}\\
1 & -\mathbf{q}^T & {\bf 0} & -z^*
\end{bmatrix} \right)
\]

Now we rewrite (1.3) using the notation from \cite{chv} page 100. Extend ${\bf c}$ to $\R^{m+n}$ by adding $m$ zeroes at the end.   Let $A_B$ be the submatrix of $A=[A_0\, I\,]$ with columns indexed by $B$, and similarly for $A_N$. Then $A_B$ is non-singular and (1.3) is of the form
\begin{equation}
\begin{aligned}
    \mathbf{x}_B &= A^{-1}_B \mathbf{b} - A^{-1}_BA_N\mathbf{x}_N\\
    z &= \mathbf{c}_B^T A^{-1}_B \mathbf{b} + (\mathbf{c}_N^T - \mathbf{c}_B^T A^{-1}_B A_N) \mathbf{x}_N
\end{aligned}
\end{equation}

It follows that $[{\bf y},w]$ is a solution to (1.4) if an only if 
\[
\begin{aligned}
[1,{\bf y}_N^T,{\bf y}_B^T,w]&\in {\rm row\ space}\left( \begin{bmatrix}
    {\bf 0}^T & A_B^{-1} A_N & I & -A_B^{-1}  \mathbf{b}\\
    1 & \mathbf{c}_B^T A_B^{-1} A_N - \mathbf{c}_N^T & {\bf 0} & -\mathbf{c}_B^T A_B^{-1} \mathbf{b}
    \end{bmatrix} \right)\\
    &={\rm row\ space}\left( \begin{bmatrix}
A_B^{-1} & {\bf 0}^T\\
\mathbf{c}_B^T A_B^{-1} & 1
\end{bmatrix}\cdot
\begin{bmatrix}
{\bf 0}^T & A_N & A_B & -\mathbf{b}\\
1 & -\mathbf{c}_N^T & -\mathbf{c}_B^T & 0
\end{bmatrix}  \right)\\
&={\rm row\ space}\left( 
\begin{bmatrix}
{\bf 0}^T & A_N & A_B & -\mathbf{b}\\
1 & -\mathbf{c}_N^T & -\mathbf{c}_B^T & 0
\end{bmatrix}  \right)\\
\end{aligned}
\]
because
$\displaystyle\begin{bmatrix}
A_B^{-1} & {\bf 0}^T\\
\mathbf{c}_B^T A_B^{-1} & 1
\end{bmatrix}
$
is nonsingular.

But this is equivalent to $[1,{\bf y}^T,w]\in{\rm row\ space}(R)$, which is what we wished to prove.


\begin{thebibliography}{10}
\bibitem{chv}
V. Chv\'{a}tal, {\em Linear Programming},  W. H. Freeman and Company, New York, 1983
\bibitem{vdb}
R. J. Vanderbei, {\em Linear Programming, Foundations and Extensions}, Springer, New York, 2014 
\bibitem{todd}
M. J. Todd, {\em Linear and Quadratic Programming in Oriented Matroids}, 
Journal of Combinatorial Theory, Series B 39, 105-133 (1985)
\end{thebibliography}
\end{document}